\documentclass[12pt]{article}
\usepackage{amssymb}
\usepackage{amsmath}
\usepackage{amsmath,amssymb,amsfonts}
\usepackage{graphicx}
\usepackage{color}

\providecommand{\U}[1]{\protect\rule{.1in}{.1in}}
\newtheorem{teo}{Theorem}

\newtheorem{lemma}{Lemma}

\newtheorem{definition}{Definition}
\newtheorem{remark}{Remark}

\begin{document}

\title{Shallow Water Model for Lakes with Navier slip boundary condition}
\author{N.V. Chemetov, F. Cipriano and S. Gavrilyuk \\
{\small DCM-FFCLRP / University of Sao Paulo, Brazil},\\
{\small GFM e Dep.\ de Matem{\'a}tica FCT-UNL, Lisboa, Portugal}, \\
{\small Marseille University and CNRS UMR 6595, IUSTI, Marseille, France}}
\date{}
\maketitle
\tableofcontents

\bigskip

\abstract{We study a  model 
describing the motion of the fluid in a lake, assuming  
 inflow-outflow effects across the bottom, the porous coast and the 
inflows and outflows of rivers.
 We prove the global in time existence result  for this model. 
The solvability is shown in the class of
solutions with  $L_{p}$-bounded vorticity for any given $p\in (1,\infty]$.}

\bigskip



\section{Modelling of the lake equations with inlets and outlets}

\label{sec0}

We consider a motion of water in a lake with slow varying bottom topography
and with inflow-outflow effects. The motion of the water is described by
Euler type equations 
\begin{equation}
\mathbf{v}_{t}+\left( \mathbf{v\cdot \nabla }\right) \mathbf{v+}g\mathbf{%
\nabla }h=-\frac{\varkappa }{b}\mathbf{v}+\frac{\mathbf{G}}{b},\quad \mathrm{%
div}\left( b\mathbf{v}\right) =A,\qquad \left( \mathbf{x,}t\right) \in
\Omega _{T}:=\Omega \times \lbrack 0,T]  \label{eq1001}
\end{equation}%
with an initial condition 
\begin{equation}
\mathbf{v}(\mathbf{x},0)=\mathbf{v}_{0}(\mathbf{x}),\quad \quad \mathbf{x}%
\in \Omega ,  \label{master2}
\end{equation}%
satisfying 
\begin{equation}
\mathrm{div}(b\,\mathbf{v}_{0})=A(\mathbf{x},0)\quad \quad \mathbf{x}\in
\,\Omega .  \label{eq1006}
\end{equation}%
Here $\Omega \subseteq \mathbb{R}^{2}$ is a bounded domain, $\mathbf{x}%
=(x,y) $ are cartesian coordinates, $\mathbf{v}=(u,v)$ is the fluid velocity
averaged over the lake depth, $g$ is the constant gravity acceleration, $h$
is the position of the free surface, $b(\mathbf{x})>0$ is the lake depth.
The coefficients $\varkappa ,\mathbf{G},A$ are known parameters,
characterizing the nature of porous properties of the bottom of the lake.
The right-hand side of the momentum equation is a friction term and $%
\varkappa (\mathbf{x},t)\geqslant 0$ is the friction coefficient. $A(\mathbf{%
x},t)$ is a so-called source term. When $\varkappa ,\mathbf{G},\,A=0$,
system (\ref{eq1001}) represents classical lake equations (see \cite{Allen}
and \cite{Camassa}). Essentially this 2-dimensional system can be derived
from the 3-dimensional Navier-Stokes equations, assuming that the ratio of
vertical to horizontal length scales is very small and also the ratio of
horizontal speeds and the gravity wave speeds is small too. The friction
term can be obtained rigorously from the Navier-Stokes equations by
averaging over the lake depth $b,$ assuming the fulfillment of a so-called
Navier slip (or friction) boundary condition on the bottom of the lake (see 
\cite{Gerbeau}, \cite{Gordon} and \cite{jager}). A real physical picture is
shown in Figure 1.

\begin{figure}[htbp]
\begin{center}
\caption[width=14cm] 
{\it An incompressible fluid is confined to a 3d basin by a 
downward gravitational field  $g$. The top level $z=h(x,y)$ of the lake is a ``rigid lid'', i.e on $ z=h $ 
we have the conditions
$\,\,
\vec{v}\cdot \vec{n}=0, \quad   2\, D(\vec{v}) \vec{n} \cdot \vec{\tau}_{i}=0.\,\,$ Here $\vec{v}\,$ 
is the 3d vector velocity of the fluid;  $\vec{\tau}_{i}$ are tangent vectors to $ z=h $ and $\vec{n}
$ is the outward normal vector. On the bottom   $z=-b(x,y)$ of the lake  we impose the conditions 
$\,\,
\vec{v}\cdot \vec{n}=A, \quad
2\mu \;\vec{\tau}_{i}\; D( \vec{v}) \vec{n}+\varkappa \,( \vec{v} \cdot \vec{\tau}_{i})=G_{i} . \,\,
$
Here $A$
are the quantity of inflow and outflow water through the bottom 
 and $\varkappa , \, G_{i}$
describe friction features of the bottom 
(roughness of the bottom: seaweed, stones). 
} \label{bfig1}
\end{center}
\includegraphics[width=15cm]{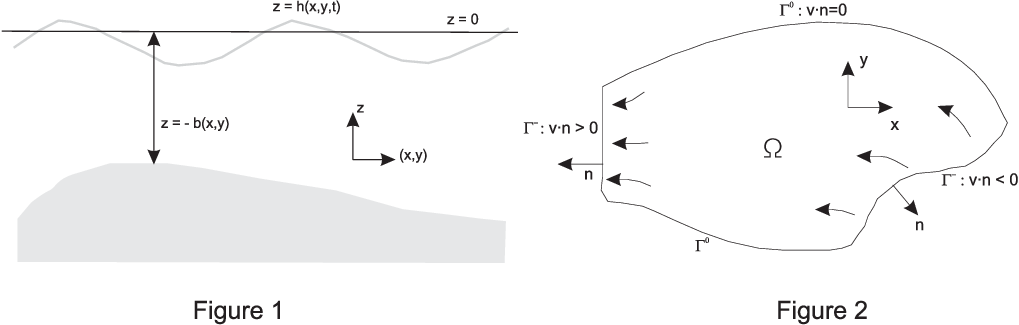}
\end{figure}

When $b=1$ and $\varkappa ,\mathbf{G},\,A=0$, these equations \eqref{eq1001}-%
\eqref{eq1006} correspond to the 2-dimensional Euler equations. The
techniques employed in the 2-dimensional Euler equations (see, for instance, 
\cite{kato}, \cite{yu}) were successfully extended to the case $b\neq 1$. The
solvability of the lake equations ($\varkappa ,\mathbf{G},\,A=0$) is very
well studied, we refer the reader to \cite{bresch}, \cite{levermore1}, \cite%
{levermore2}, \cite{oliver}. In these mentioned articles the non-penetration
condition $\mathbf{v}\cdot \mathsf{n}=0$ was assumed on the boundary $\Gamma 
$ of the domain $\Omega $. Here $\mathsf{n}$ is the normal vector to $\Gamma 
$. In the present article we investigate a more general case, when $%
\varkappa ,\mathbf{G},\,A$ are non zero terms and there are inlets and
outlets of the fluid through $\Gamma $, i.e. the lake has a porous coast,
inflows and outflows of rivers (see Figure 2).

Let the fluid flow into and out the domain $\Omega $ through parts $\Gamma
^{-}$ and $\Gamma ^{+}$ of the boundary $\Gamma $, such that $\Gamma =\Gamma
^{-}\cup \Gamma ^{0}\cup \Gamma ^{+}$ and $meas(\Gamma ^{+})\neq
0,\,meas(\Gamma ^{-})\neq 0$. Let the quantity of the inflowed and outflowed
fluid be equal to $a$ 
\begin{equation}
\mathbf{v}\cdot \mathsf{n}=a\quad \quad \mbox{ on }\Gamma _{T},
\label{eq1.39}
\end{equation}%
with 
\begin{equation}
a(\mathbf{x},t)\ \left\{ \begin{aligned} < 0 ,\quad \quad & \quad\text{ if
}(\mathbf{x},t)\in \Gamma_{T}^{-}:=\Gamma^{-}\times (0,T), \\ = 0 ,\quad
\quad & \quad\text{ if }(\mathbf{x},t)\in \Gamma_{T}^{0}:=\Gamma^{0}\times
(0,T), \\ > 0 ,\quad \quad & \quad\text{ if }(\mathbf{x},t)\in
\Gamma_{T}^{+}:=\Gamma^{+}\times (0,T) , \end{aligned}\right.
\label{apositive}
\end{equation}%
satisfying 
\begin{equation}
\int_{\Gamma }b(\mathbf{x})\,a(\mathbf{x},t)\,\,d\mathbf{x}=\int_{\Omega }A(%
\mathbf{x},t)\,\,d\mathbf{x}\quad \quad \mbox{ for }\,t\in \lbrack 0,T].
\label{eqC2}
\end{equation}%
Let us rewrite the system \eqref{eq1001}-\eqref{eq1006} in terms of the
velocity $\mathbf{v}$ and the vorticity $\mathrm{rot}(\mathbf{v}%
):=v_{x}-u_{y}$. Let us define by $\nabla ^{\bot }:=(-\partial _{y},\partial
_{x})$ and $\omega :=\frac{1}{b}\,\mathrm{rot}(\mathbf{v})$. Using the
identity $\;\;\mathrm{rot}\left( (\mathbf{v}\cdot \nabla )\mathbf{v}\right) =%
\mathrm{div}(\mathbf{v}\,\mathrm{rot}(\mathbf{v}))$, we see that our system
can be written as 
\begin{align}
\partial _{t}(b\omega )+\mathrm{div}(b\,\omega \,\mathbf{v})& =-\varkappa
\omega -(\mathbf{v}\cdot \nabla ^{\perp })\left( \frac{\varkappa }{b}\right)
+\mathrm{rot}\bigl(\frac{\mathbf{G}}{b}\bigl)\quad \quad \text{in}\,\Omega
_{T},  \label{eq888} \\
\mathrm{rot}(\mathbf{v})& =b\,\omega ,\quad \quad \quad \quad \mathrm{div}%
(b\,\mathbf{v})=A\quad \quad \quad \quad \quad \text{in}\,\Omega _{T}
\label{eq1010}
\end{align}%
with the initial condition 
\begin{equation}
\omega (\mathbf{x},0)=\omega _{0}(\mathbf{x}):=\frac{1}{b}\mathrm{rot\ }%
\mathbf{v}_{0}(\mathbf{x}),\quad \quad \mathbf{x}\in \Omega .  \label{7eq7}
\end{equation}%
On the part $\Gamma _{T}^{-}$ we impose the Navier slip boundary condition 
\begin{equation}
2D(\mathbf{v})\mathsf{n}\cdot \mathsf{s}+\alpha \mathbf{v}\cdot {\mathsf{s}}%
=\eta \quad \quad \mbox{ on }\Gamma _{T}^{-},  \label{eqq1.2}
\end{equation}%
where $\alpha ,\eta $ are known functions, describing the porous properties
of $\Gamma ^{-}$ and $\mathsf{s}$ is the tangent vector to $\Gamma $. The
tensor $D(\mathbf{v}):=\frac{1}{2}[\nabla \mathbf{v}+(\nabla \mathbf{v}%
)^{T}] $ is the rate-of-strain.

In the article \cite{C&A} the well posedness for Euler equations
with the Navier slip boundary conditions \eqref{eq1.39} and \eqref{eqq1.2}
for $L_{p}-$ bounded vorticity when $p>2$ was proved.   In the present article we apply
the method developed in \cite{C&A} to prove the solvability of the problem 
\eqref{eq1001}-\eqref{eq1006}, \eqref{eq1.39}-\eqref{eqC2}, \eqref{eqq1.2},
extending on the case of $L_{p}-$bounded vorticity when $p\in (1,\infty ].$  We refer the reader for the hystory of the 
solvability of the Euler equations to \cite{C&A}. The actual article is a version of the article which was published in \cite{CC6}.
   Since the publication of the papers \cite{C&A} and \cite{CC6}, other articles on the subject studied in these articles have been written, such as \cite{CCC3}-\cite{CC5}, where the boundary layer theory have been developed for Navier-Stokes and Euler equations with Navier slip boundary condition. We also refer to the problem for Euler equations with sources and sinks, which was considered by Chemetov, Starovoitov \cite{chem}.  

A significant progress in the study of the conservative non-linear hyperbolic law coupled with elliptic equations, being similar systems as appeared in the formulation of  the Euler equations in terms of vorticity - stream function. We refer to the articles \cite{ant2, AC2},  \cite{chem1, chem2}, \cite{chem3}-\cite{chem5}, where the Kruzkov approach  has been generalized  using the kinetic method.  The stochastic non-linear hyperbolic-elliptic systems have been investigated in \cite{ACC21}, \cite{CC0}.

\section{The existence result}

\setcounter{equation}{0}

\label{sec444}
We shall introduce some notations from \cite{LadySolonUral68}, \cite%
{LadyUral68}. Let $l$, $m$ be non-negative reals and $q\geq 1$; we 
consider  the Sobolev spaces $L_{q}(\Omega )$, $%
W_{q}^{l}(\Omega )$, $W_{q}^{l,\,m}(\Omega_{T})$, $W_{q}^{l}(\Gamma )$, $%
W_{q}^{l,\,m}(\Gamma_{T})$ and the H\"{o}lder spaces
$C^{l,\,m}(\Omega_{T})$, $C^{l,\,m}(\Gamma_{T})$, where $l$ and $m$
 correspond to the regularity on the variables $\mathbf{x}$
and $t$, respectively.

We assume that the boundary $\Gamma $ is $C^{2}$-smooth. If we parametrize
the boundary $\Gamma $ using the arc length $s$, it follows that $\frac{d%
\mathsf{n}}{ds}=k{\mathsf{s}}$, where the curvature $k$ of $\Gamma $ is
continuous function on $\Gamma .$

We consider the following  characterization of the Navier slip boundary condition
\eqref{eqq1.2}.
\begin{lemma}
\label{lem4.1} Let $V_{a}:=\{\mathbf{w}\in C^{2}(\bar{\Omega}):\,\mathbf{w}%
\cdot \mathsf{n}=a\mbox{ on
}\Gamma \}.$ A vector $\mathbf{u}\in V_{a}$ satisfies the condition %
\eqref{eqq1.2} if and only if $\quad \quad \mathrm{rot}(\mathbf{u}%
)=(2k-\alpha )\,\mathbf{u}\cdot \mathsf{s}+(\eta -2\,a_{s}^{\prime }).$
\end{lemma}

\vspace{1pt}

This lemma is a generalization of Lemma 4.1 and Corollaries 4.2, 4.3 of \cite{kel}. 
Let us introduce the following functions on the boundary $\Gamma _{T}$
\begin{equation*}
\gamma :=\frac{1}{b}\bigl(2k-\alpha \bigr),\quad g:=\frac{1}{b}\bigl(\eta
-2\,a_{s}^{\prime }\bigr),\quad \omega _{\Gamma }(\mathbf{u}):=\gamma \,%
\mathbf{u}\cdot {\mathsf{s}}+g.
\end{equation*}%
We can observe that the problem for the velocity given by 
 the system \eqref{eq1001}-\eqref{eqC2}, \eqref{eqq1.2} is
equivalent to the problem for the velocity $\mathbf{v}$ and the modified
vorticity $\omega $, given by  the system \eqref{eq888} -\eqref{eq1010}
with the initial condition \eqref{7eq7} and the boundary condition 
\begin{equation}
\omega =\omega _{\Gamma }(\mathbf{v})\quad \text{on}\quad \Gamma _{T}^{-}.
\label{eqcf1}
\end{equation}

Let us suppose just for a moment that $\omega $ is a known function, then
the solution $\mathbf{v}$ of the system \eqref{eq1.39}, \eqref{eq1010} can
be written in the form 
\begin{equation}
\mathbf{v}=\frac{1}{b}\nabla ^{\perp }h+\nabla H,\quad \quad \quad (\mathbf{x%
},t)\in \Omega _{T},  \label{V}
\end{equation}%
where $h$, $H$ are solutions of the following systems 
\begin{equation*}
\left\{ 
\begin{array}{ll}
-\text{div}(\frac{1}{b}\nabla h)=b\,\omega , & (\mathbf{x},t)\in \Omega _{T},
\\ 
&  \\ 
h=0, & (\mathbf{x},t)\in \Gamma _{T},%
\end{array}%
\right. \,\mbox{ }\,\left\{ 
\begin{array}{ll}
\text{div}(b\nabla H)=A, & (\mathbf{x},t)\in \Omega _{T}, \\ 
&  \\ 
\frac{\partial H}{\partial n}=a, & (\mathbf{x},t)\in \Gamma _{T}.%
\end{array}%
\right.
\end{equation*}
From classical results of potential theory for elliptic equations (see \cite%
{vlad}), the function $h$ can be represented in the integral form 
\begin{align}
h(\mathbf{x},t) =\int_{\Omega }K_{1}(\mathbf{x},\mathbf{y})\,\,b(\mathbf{y}%
)\,\omega (\mathbf{y},t)\,d\mathbf{y}=:K_{1}\ast b\,\omega  \label{eq66sec2}
\end{align}
and from the theory for elliptic equations (see \cite{LadyUral68}) we have

\begin{lemma}
\label{teo2sec2} For any given $q\in (1,\infty)$ and any fixed $t\in[0,T]$,
the functions $h$ and $H$ satisfy the following estimates 
\begin{align}
||h(\cdot,t)||_{W_{q}^{2}(\Omega )}& \leqslant
C||b\,\omega(\cdot,t)||_{L_{q}(\Omega)},  \label{ell-1} \\
||H(\cdot,t)||_{W_{q}^{2}(\Omega )}& \leqslant C (
\|A(\cdot,t)\|_{L_q(\Omega)}+\|a(\cdot,t)\|_{W_{q}^{1-\frac{1}{q}}(\Gamma )}
).  \label{ell-6}
\end{align}
\end{lemma}

In view of \eqref{V}, the function $\mathbf{v}$ fulfills 
\begin{equation}  \label{VELOC}
\|\mathbf{v}(\cdot,t )\|_{W^1_{q}(\Omega)} \leqslant C ( \| \omega (\cdot,
t)\|_{L_q(\Omega)} +\|A(\cdot, t)\|_{L_{q}(\Omega)}+ \|a(\cdot, t
)\|_{W_{q}^{1-\frac{1}{q} }(\Gamma)} )
\end{equation}
for any given $q\in (1,\infty)$ and any $t\in[0,T]$.

\bigskip

We consider $p\in (1,\infty ]$  and suppose that the  data satisfy the
following conditions 
\begin{equation}
\begin{cases}
0<b\in C^{1}(\bar{\Omega}),\hspace{3cm} &  \omega _{0}\in L_{p}(\Omega ) ,\\ 
A\in
L_{2 }(0,T,W_{\widetilde{p}}^{1}(\Omega )),  & a\in
L_{2 }(0,T,W_{p}^1 (\Gamma )) ,\\
\mathrm{rot}(\mathbf{G})\in L_{1 }(0,T,L_{p}(\Omega )), & \varkappa \in L_{\infty }(\Omega _{T}), \\ 
\alpha \in L_{2 }(0,T,L_{\widetilde{p}}(\Gamma^- )), & \eta \in L_{1
}(0,T,L_{p}(\Gamma ^{-})) 
\end{cases}
\label{eq00sec1}
\end{equation}%
with 
\begin{equation}
\widetilde{p}=%
\begin{cases}
p,\quad  & \text{if $p>2$}, \\ 
2+\varepsilon \quad \,\,\text{with some $\varepsilon >0$,} & \text{if $p=2$,}
\\ 
\frac{p}{p-1}, & \text{if $1<p<2$.}%
\end{cases}
\label{eqlinha}
\end{equation}

\bigskip

Next, we follow \cite{Delort}, \cite{schochet} in order to introduce the  notion 
of a weak solution  of our problem.

\begin{definition}
\label{def2sec1} A pair of functions $\{\omega ,\mathbf{v}\}$ is said to be
a \underline{weak solution} of the problem \eqref{eq1001}-\eqref{eqC2}, %
\eqref{eqq1.2}, if $\omega \in L_{\infty }(0,T,L_{p}(\Omega )),\,$ $\mathbf{v%
}\in L_{\infty }(0,T,\,W_{q}^{1}(\Omega ))$ (here $q=p,$\ if $p\in (1,\infty
),$ any $q\in \lbrack 1,\infty )$ if $p=\infty )$ and the following
equalities hold 
\begin{equation}  \label{eq8}
\mathrm{rot}(\mathbf{v})=b\,\omega , \; \quad \mathrm{div}(b\,\mathbf{v})=A
\quad \text{ a.e. in }\,\Omega _{T} \quad \text{ and } \quad \mathbf{v}\cdot 
\mathsf{n}=a \quad \text{ a.e. on }\,\Gamma _{T},
\end{equation}
\begin{align}  \label{eq77}
& \int_{\Omega _{T}}b\,\omega (\psi _{t}+\nabla H\cdot \nabla \psi )\,d%
\mathbf{x}dt+F[\omega ,K_{\psi }]-\bigl[\varkappa \omega +(\mathbf{v}\cdot
\nabla ^{\perp })\bigl(\frac{A}{b}\bigr)-\mathrm{rot}\bigl(\frac{\mathbf{G}}{%
b}\bigr)\bigr]\psi \,d\mathbf{x}dt  \notag \\
& =-\int_{\Omega }b\,\omega _{0}\,\psi (\mathbf{x},0)\,d\mathbf{x}%
+\int_{\Gamma _{T}^{-}}b\,a\,\omega _{\Gamma }(\mathbf{v})\,\psi \,d\mathbf{x%
}dt.
\end{align}
where 
\begin{equation}
F[\omega ,K_{\psi }]:=\int_{0}^{T}\int_{\Omega }\int_{\Omega }\omega (%
\mathbf{x},t)\omega (\mathbf{y},t) \; K_{\psi }(\mathbf{x},\mathbf{y},t)\,d%
\mathbf{x}d\mathbf{y}dt  \label{WW}
\end{equation}%
and 
\begin{equation*}
K_{\psi }(\mathbf{x},\mathbf{y},t):=\nabla _{\mathbf{x}}^{\perp }K_{1}(%
\mathbf{x},\mathbf{y})\; \frac{b(\mathbf{y})\nabla _{\mathbf{x}}\psi (%
\mathbf{x},t)-b(\mathbf{x})\nabla _{\mathbf{y}}\psi (\mathbf{y},t)}{2}
\end{equation*}%
for an arbitrary function $\psi \in C^{1,1}(\overline{\Omega }_{T}),$
satisfying the condition 
\begin{equation}
\mbox{ supp}\,(\psi )\subset (\Omega \cup \Gamma ^{-})\times \left[
0,T\right) .  \label{psi1}
\end{equation}%
The function $\psi $ will be called a \underline{test function}.
\end{definition}

We  remark  that this definition is an extention of the classical one; in
fact, for $\frac{4}{3}\leqslant p$,  equality \eqref{eq77} is equivalent
to the identity 
\begin{equation}
\begin{split}
& \int_{\Omega _{T}}b\,\omega (\psi _{t}+\mathbf{v}\cdot \nabla \psi )-\bigl[%
\varkappa \omega +(\mathbf{v}\cdot \nabla ^{\perp })\bigl(\frac{A}{b}\bigr)-%
\mathrm{rot}\bigl(\frac{\mathbf{G}}{b}\bigr)\bigr]\psi \,d\mathbf{x}dt \\
& =-\int_{\Omega }b\,\omega _{0}\,\psi (\mathbf{x},0)\,d\mathbf{x}%
+\int_{\Gamma _{T}^{-}}b\,a\,\omega _{\Gamma }(\mathbf{v})\,\psi \,d\mathbf{x%
}dt.
\end{split}
\label{eq7}
\end{equation}%
We establish our main result.

\begin{teo}
\label{teo2sec1} Suppose that the data $b,\,\,a,\,\,\alpha ,\,\,\eta ,\,\,\omega
_{0},\,\,A,\,\mathbf{G}$ satisfy \eqref{eq00sec1}-\eqref{eqlinha}, then
there exists at least one weak solution $\{\omega ,\mathbf{v}\}$ of the
problem \eqref{eq1001}-\eqref{eqC2}, \eqref{eqq1.2}.
\end{teo}

\begin{remark}
This existence result is based on  inequality \eqref{A}, which can be
 obtained from \eqref{eq888}, if we assume that $\,\omega $ is a regular
function. But in reality $\omega \in L_{\infty }(0,T,L_{p}(\Omega ))$, so the  trace
on $\Gamma _{T}^{-}$ is not defined.
 Fortunately $\omega $ satisfies the
equation \eqref{eq888}, which can be represented in the form $div_{\mathbf{x}%
,t}(\mathbf{F})=z,$ allowing to  define in a correct way  the trace of $\omega 
$\ \ on $\Gamma _{T}^{-}$ and for $t=0,$ according to the formulas %
\eqref{bo} and \eqref{to}. For the reader we refer the book of Temam \cite%
{Kuf}, pages 6-9, where it is perfectly described this trace approach. Also in
the survey paper \cite{chen} the theory of trace is discussed for divergence
measure fields.
\end{remark}

For the reader's convenience, we recall some well-known  embedding results.
 By Sobolev's embedding theorems
p. 287 \cite{Kuf}, we have

\begin{lemma}
Let $Q\subseteq \mathbb{R}^{n},\quad n\geqslant 1$ be a locally lipschitzian
domain. For any $q\in \lbrack 1,\infty )$, the following continuous
embedding holds 
\begin{equation}
W_{q}^{1}(Q)\hookrightarrow \left\{ \begin{aligned} C^{\alpha }( \bar Q )
,\quad &\text{ if } \quad q > n \quad \,\, \text{ for } \quad \quad \,\,
\alpha = 1-\frac{n}{q} , \\ L_{\widetilde q}(Q ) ,\quad &\text{ if } \quad
q=n \quad \:\: \text{ for any } \quad \widetilde{q}\in [1, \infty ), \\
L_{\widetilde q}(Q ) ,\quad &\text{ if } \quad q < n \,\, \quad \text{ for
any }\quad \widetilde{q}\in [1, \frac{ n q}{n- q } ]. \end{aligned}\right.
\label{ell-90}
\end{equation}
\end{lemma}

For traces on the boundary $\Gamma $ of our domain $\Omega $
( see theorems p. 287, p. 319 and p. 336 \cite{Kuf}),
 we have 

\begin{lemma}
For any $q\in [1,\infty)$ 
\begin{equation}
W_{q}^{1}(\Omega )\hookrightarrow \left\{ \begin{aligned}
C^\alpha(\bar{\Omega}),\quad &\text{ if } \quad q> 2 \quad \quad \quad
\text{ for }\,\, \quad \quad \alpha=1-\frac{2}{q},\\
L_{\widetilde{q}}(\Gamma),\quad &\text{ if } \quad q=2\quad \quad \quad
\text{ for any } \quad\widetilde{q}\in [1, \infty ),\\
L_{\widetilde{q}}(\Gamma),\quad &\text{ if } \quad q < 2 \quad \quad \quad
\text{ for any }\quad \widetilde{q} \in [1, \frac{q}{2-q} ] . \end{aligned} %
\right.  \label{ell-9}
\end{equation}
\end{lemma}

\section{Solvability of Navier-Stokes type system}

\setcounter{equation}{0}

To prove  Theorem 1 we follow the methods in \cite{C&A}.
A solution will be obtained by the limit  of viscous solutions
 of a Navier-Stokes
type system with artificial viscosity $\nu>0$.
By this reason,  we sketch  the proof
of the solvability of such corresponding viscous system. 
 Let us remember a
result from  approximation theory.

\begin{lemma}
\label{ap1} Let $Q\subseteq \mathbb{R}^{n},\quad n\geqslant 1$ be an open
set and $q\in \lbrack \,1,+\infty ]$. Then for any $b\in L_{q}\,(Q)$ there
exist functions $b^{\theta }\in C^{\infty }(Q),$ satisfying the following
properties: 
\begin{equation}
||b^{\theta }||_{L_{q}(Q)}\leqslant C||b||_{L_{q}(Q)}\,,  \label{ap2}
\end{equation}%
\begin{equation}
b^{\theta }\underset{{\theta }\rightarrow 0}{\longrightarrow }b\quad \quad %
\mbox{ in }\quad \quad L_{q^{\prime }}(Q),  \label{ap3}
\end{equation}%
for $q^{\prime }=q,$ if $q<\infty $ and for any $q^{\prime }<\infty $, if $%
q=\infty .$
\end{lemma}

Using Lemma \ref{ap1}, we approximate the data $a,\,\gamma ,\,g,\,\omega
_{0},\,\varkappa ,\,A,\,\mathbf{G}$ by $C^{\infty }$- 
functions $a^{\theta },\,\gamma ^{\theta },\,g^{\theta },\,\omega
_{0}^{\theta },\,\varkappa ^{\theta },\,A^{\theta },\,\mathbf{G}^{\theta }$,
according to the relations \eqref{ap2}-\eqref{ap3} in the respective space $%
L_{q}(Q)$, defined by the regularity conditions \eqref{eq00sec1}-\eqref{eqlinha}.
We can assume that any derivatives of these smooth approximations are bounded by
constants, just  depending  on a \underline{\textit{fixed}} parameter $\theta $
and that $a^{\theta }$ fulfills the relation \eqref{apositive}. The
functions $\gamma ^{\theta },\,\,g^{\theta }$, $\omega _{0}^{\theta }$ satisfy
a so-called \textit{compatibility} condition 
\begin{equation*}
\gamma ^{\theta }(\cdot ,t),\,\,g^{\theta }(\cdot ,t)=0\quad \mbox{ for }t\in %
\left[ 0,\theta \right] \quad \quad \mbox{ and }\quad \quad \omega
_{0}^{\theta }(\mathbf{x})=0\quad \mbox{ for }\mathbf{x}\in U_{\theta
}(\Gamma ).
\end{equation*}%
In the following considerations we shall work with these approximations of the  data,
but for  sake of simplicity in the notations, we suppress the dependence
on ${\theta }$ and continue to write $a,\,\gamma ,\,g,\,\omega _{0},\,A,\,%
\mathbf{G}$, respectively.

Next, we construct the pair $\{\omega ,h\}$ as a solution of an auxiliary 
\textbf{Problem}, which is a coupling of two following systems. \bigskip

\textbf{Problem}. \textit{Find} $\omega \in W_{2}^{2,1}(\Omega _{T})$, 
\textit{satisfying the system} 
\begin{equation}
\begin{cases}
\displaystyle{\ \partial _{t}(b\omega )+\mathrm{div}(b\,\omega \,\mathbf{v}%
)=-\varkappa \omega -(\mathbf{v}\cdot \nabla ^{\perp })\bigl(\frac{A}{b}%
\bigr)+\mathrm{rot}\bigl(\frac{\mathbf{G}}{b}\bigr)}+\nu \Delta \omega \quad
\quad \text{on $\Omega _{T}$,} \\ 
\omega \big|_{\Gamma _{T}}=\omega _{\Gamma }(\mathbf{v}),\quad \quad \quad
\quad \omega \big|_{t=0}=\omega _{0}%
\end{cases}
\label{eq3sec2}
\end{equation}%
\textit{and find} $h,H\in L_{\infty }(0,T,\ W_{2}^{2}(\Omega )))$, \textit{%
satisfying the systems } 
\begin{equation}
\left\{ 
\begin{array}{ll}
-\mathrm{div}\bigl(\frac{1}{b}\nabla h\bigr)=b\langle \omega \rangle , & (%
\mathbf{x},t)\in \Omega _{T}, \\ 
&  \\ 
h=0, & (\mathbf{x},t)\in \Gamma _{T},%
\end{array}%
\right. \,\mbox{ }\,\left\{ 
\begin{array}{ll}
\mathrm{div}\bigl(b\,\nabla H\bigr)=A, & (\mathbf{x},t)\in \Omega _{T}, \\ 
&  \\ 
\frac{\partial H}{\partial n}=a, & (\mathbf{x},t)\in \Gamma _{T}%
\end{array}%
\right.  \label{eq4sec2}
\end{equation}%
\textit{and the velocity }${\,\mathbf{v}\,}$\textit{\ satisfies the relation %
\eqref{V}. Here we denoted by} 
\begin{equation}
\left\langle \omega \right\rangle (\cdot ,t):=\frac{1}{\theta }%
\int_{t-\theta }^{t}\mathit{\ }\left[ \omega (\cdot ,s)\right] _{R}\ ds
\label{R}
\end{equation}%
\textit{with the cut-off function} $\left[ \cdot \right] _{R}$, \textit{%
defined as} $\left[ \omega \right] _{R}:=\max \Big\{-R,\,\min \{R,\omega \}%
\Big\}$. \textit{We assume that} $\omega =0$ \textit{outside of the interval}
$\left[ 0,T\right] $. \bigskip

The solution $\{\omega ,h\}$ of \textbf{Problem} depends on the parameters $%
\nu ,\ \theta$, but for the simplicity of presentation  in the sequel we just 
indicate the dependence of functions and constants on the parameters $\nu ,
\ \theta $, if it is necessary.

Using Schauder's fixed point argument developed in \cite{C&A}, we can show
the following Lemma.

\begin{lemma}
\label{lem3sec2} There exists at least one \underline{weak} solution $%
\{\omega ,h\}$ of the systems (\ref{eq3sec2})-(\ref{eq4sec2}), such that for
some $\alpha \in (0,1)$ 
\begin{equation}
\omega \in V_{2}^{1,0}(\Omega _{T})\cap C(0,T,\,L_{2}(\Omega )),\quad \quad
h\in C(0,T,\,C^{1+\alpha }(\overline{\Omega })).  \label{eq3.6}
\end{equation}
\end{lemma}

In this lemma the space $V_{2}^{1,0}(\Omega _{T})$ has the norm $%
||u||_{V_{2}^{1,0}(\Omega _{T})}=\max_{t\in \lbrack 0,T]}||u(\mathbf{x}%
,t)||_{L_{2}(\Omega )}+||u_{\mathbf{x}}(\mathbf{x},t)||_{L_{2}(\Omega
_{T})}. $

\bigskip

Now we shall deduce  a priori estimates for the solution $\{\omega ,h\}$ 
of \textbf{%
Problem} $\mathbf{P},$ independent of  ${\nu \in (0,1)}$, therefore
it is convenient to write $\{{\ \omega_{\nu},h_{\nu }}\}$
 and $\mathbf{v_{\nu}}=\frac{1%
}{b}\nabla^\perp h_{\nu}+\nabla H$, the  parameter $\theta$ being fixed.

\begin{teo}
\label{teo4sec2} There exists at least one weak solution $\{\omega _{{\nu }%
},h_{\nu }\}$ of \textbf{Problem} $\mathbf{P}$, which satisfies the
following 
\begin{eqnarray*}
\omega _{{\nu }} &\in &W_{2}^{2,1}(\Omega _{T})\cap C^{\alpha ,\alpha /2}(%
\overline{\Omega }_{T}), \\
h_{\nu } &\in &L_{\infty }(0,T,\ W_{q}^{2}(\Omega ))\cap C^{2+\alpha ,\alpha
/2}(\overline{\Omega }_{T}),\quad \quad \partial _{t}\,h_{\nu }\in L_{\infty
}(0,T,\ W_{q}^{2}(\Omega ))
\end{eqnarray*}%
for any fixed $q\in (1,\infty )$ and for some $\alpha \in (0,1)$ and 
\begin{align}
||\omega _{{\nu }}||_{L_{\infty }(\Omega _{T})}& \leqslant B,  \label{w2} \\
||h_{\nu }||_{L_{\infty }(0,T,\,W_{q}^{2}(\Omega ))}& \leqslant C,\quad
||h_{\nu }||_{L_{\infty }(0,T,\,C^{1+\alpha }(\overline{\Omega }))}\leqslant
C,  \label{h1} \\
||\partial _{t}h_{\nu }||_{L_{\infty }(0,T,\,W_{q}^{2}(\Omega ))}& \leqslant
C,\quad ||\partial _{t}h_{\nu }||_{L_{\infty }(0,T,\,C^{1+\alpha }(\overline{%
\Omega }))}\leqslant C  \label{h2}
\end{align}%
with the constants $B=B(\theta )$, $C=C(\theta ,q),$ which are independent
of $\nu $.
\end{teo}

\noindent

\textbf{Proof.} Since $b\,\langle\omega_\nu\rangle$ and $\frac{\partial}{%
\partial t} \bigl(b\,\langle\omega_\nu\rangle\bigr)$ are bounded functions
by a constant $C=C(\theta)$ in the $L_\infty$-norm, then  estimates %
\eqref{h1} and \eqref{h2} are a direct consequence of \eqref{eq4sec2}, %
\eqref{ell-1} and \eqref{ell-90}.

On the other hand, the function $\omega _{{\nu }}$ is a solution of the equation 
\begin{equation*}
\partial _{t}\omega _{{\nu }}-{\nu }\Delta \omega _{{\nu }}=F(\mathbf{x}%
,t),\quad \quad (\mathbf{x},t)\in \Omega _{T}
\end{equation*}%
with $F:=-(\mathbf{v}_{\nu }\cdot \nabla )\omega _{\nu }-\frac{A+\varkappa }{%
b}\omega _{\nu }-\frac{1}{b}(\mathbf{v}_{\nu }\cdot \nabla ^{\perp })\left( 
\frac{A}{b}\right) +\frac{1}{b}\mathrm{rot}\left( \frac{\mathbf{G}}{b}%
\right) \in L_{2}(\Omega _{T})$ by the second system of (\ref{eq4sec2}), (%
\ref{h1}). Applying Theorem 6.1, p.178, \cite{LadySolonUral68} we have $%
\omega _{{\nu }}(\mathbf{x},t)\in W_{2}^{2,1}(\Omega _{T})$, hence $\omega
_{\nu }$ satisfies equation of \eqref{eq3sec2} for a.e. $(\mathbf{x}%
,t)\in \Omega _{T}$.

To deduce (\ref{w2}) we use the maximum principle. Let

\begin{equation*}
k(t):=\max \left[ \ \max_{\overline{\Omega }}\ |\omega _{0}(x)|,\quad
||\gamma ||_{L_{\infty }(\Gamma _{T})}\cdot \max_{\Gamma }|\mathbf{v}_{{\nu }%
}(\cdot ,t)|+||g||_{L_{\infty }(\Gamma _{T})}\right] .
\end{equation*}%
Taking into account \eqref{R}, \eqref{VELOC} and \eqref{ell-90}, applied for $%
Q:=\Omega $ and any $q>2$, the velocity $\mathbf{v}_{{\nu }}$ satisfies 
\begin{equation*}
\begin{split}
\Vert \mathbf{v}_{{\nu }}(\cdot ,t)\Vert _{C(\bar{\Omega})}& \leq C(\Vert
\langle \omega _{{\nu }}(\cdot ,t)\rangle \Vert _{L_{q}(\Omega )}+1) \\
& \leq C\bigl(\int_{0}^{t}\Vert \omega _{{\nu }}(\cdot ,s)\Vert _{L_{\infty
}(\Omega )}ds+1\bigr)=:C\,\mathcal{A}(t)\quad \text{for}\quad \forall t\in
\lbrack 0,T],
\end{split}%
\end{equation*}%
where the constants $C=C(\theta )$ do not depend on $\nu $. We see that for $%
t\in \lbrack 0,T]$ 
\begin{equation}
k(t)\,,\,\left\Vert \varkappa +A+(\mathbf{v}_{\nu }\cdot \nabla ^{\perp
})\left( \frac{A}{b}\right) -\mathrm{rot}\left( \frac{\mathbf{G}}{b}\right)
\right\Vert _{L_{\infty }(\Omega )}(t)\leq C\,\mathcal{A}(t).  \label{eqestA}
\end{equation}%
Therefore multiplying the equation of \eqref{eq3sec2} by 
\begin{equation*}
q\,[\omega -k]_{+}^{q-1}:=\left\{ \begin{aligned}
q\,(\omega-k)^{q-1},\quad&\text{if}\quad\omega>k,\\
0,\quad&\text{if}\quad\omega\leqslant k, \end{aligned}\right.
\end{equation*}%
for any $q>2$, using the first equation of \eqref{eq8} and integrating over $%
\Omega $, we deduce for $y(t):=\bigl(\int_{\Omega }b\,[\omega -k]_{+}^{q}\,dx%
\bigr)^{\frac{1}{q}}\,\,$ the following inequality 
\begin{equation}
\frac{d}{dt}y^{q}(t)\leqslant q\,C\,y^{q}(t)+C\,\mathcal{A}^{q}(t)\quad 
\text{for}\quad \forall t\in \lbrack 0,T].  \label{yyyy}
\end{equation}%
Here we used the inequality 
\begin{equation}
ab\leqslant \frac{a^{r}}{r}+\frac{b^{s}}{s}\quad \text{with}\quad r,s>1,\,%
\frac{1}{r}+\frac{1}{s}=1,  \label{eqest}
\end{equation}%
taking $a:=\mathcal{A}(t)$, $b:=[\omega -k]^{q-1}$ and $r:=q$, $s:=\frac{q}{%
q-1}$. Integrating the inequality \eqref{yyyy} and taking $q\rightarrow
+\infty $, we obtain 
\begin{equation*}
\omega _{\nu }(\mathbf{x},t)\leqslant k(t)+C\,\mathcal{A}(t)\leq C\,\mathcal{%
A}(t),\quad \forall (\mathbf{x},t)\in \Omega _{T}.
\end{equation*}%
In the  same way we obtain $-C\,\mathcal{A}(t)\leqslant \omega _{\nu }(%
\mathbf{x},t),\quad \forall (\mathbf{x},t)\in \Omega _{T},$ and applying 
Gronwall's inequality to $\Vert \omega _{\nu }(\cdot ,t)\Vert _{L_{\infty
}(\Omega )}$ we deduce \eqref{w2}.

Finally, using Theorem 10.1, p.204, \cite{LadySolonUral68}, we have $%
\,\omega _{{\nu }}\in C^{\alpha ,\alpha /2}(\overline{\Omega }_{T})$ for
some $\alpha \in (0,1).$ Moreover by the theory of elliptic equations \cite%
{LadyUral68}, we conclude that $\,\,h_{\nu }\in C^{2+\alpha ,\alpha /2}(%
\overline{\Omega }_{T}).$ \bigskip $\hfill \;\blacksquare $

Choosing $R:=B$, we can take off the subscript $R$ in \eqref{R} and in the
sequel we consider that 
\begin{equation*}
\langle \omega\rangle(\cdot, t)=\frac{1}{\theta}\int_{t-\theta}^t
\omega(\cdot, s) \,ds.
\end{equation*}

\bigskip

\section{Limit transition on the viscosity}

\label{sec45} \setcounter{equation}{0}

\bigskip

In this section we consider the viscosity solutions constructed in the  previous one 
and prove the existence of the limit when  $\nu \rightarrow 0$; the
parameter $\theta$ continues to be fixed. We also prove a Gronwall type inequaly that will be 
very usefull in next section. 

\begin{lemma}
\label{lem6sec42340} For any fixed $\theta >0$ there exists a solution $%
\{\omega ,h,\mathbf{v}\}$ for a couple systems (\ref{eq7}), (\ref{eq4sec2}),
depending on ${\theta ,}$ satisfying the estimates (\ref{w2})-(\ref{h2}). Moreover, 
\textit{the velocity }${\,\mathbf{v}\,}$\textit{\ satisfies the relation %
\eqref{V}}.
\end{lemma}

\textbf{Proof.} The limit transition on $\nu \rightarrow 0$ will be done in
two steps:

\underline{1st step.} \vspace{1pt}From (\ref{w2})-(\ref{h2}) we conclude
that there exists a subsequence of $\{\omega _{{\nu }},h_{\nu }\}$, such
that 
\begin{align}
h_{\nu }\rightharpoonup h,\quad & \partial _{t}h_{\nu }\rightharpoonup
\partial _{t}h\quad & & \text{weakly}-\ast \text{ in }L_{\infty
}(0,T,\,W_{q}^{2}(\Omega )),  \label{eq31sec3} \\
\omega _{{\nu }}\rightharpoonup \omega ,\quad & \langle \omega _{\nu
}\rangle \rightharpoonup \langle \omega \rangle \quad & & \text{weakly}-\ast 
\text{ in }L_{\infty }(\Omega _{T}),  \notag
\end{align}%
the limit functions $\{\omega ,h\}$ fulfill the estimates (\ref{w2})-(\ref%
{h2}). This also implies (\ref{V}). By (\ref{eq31sec3}) and the
representation (\ref{eq66sec2}), we have 
\begin{equation}
h=K_{1}\ast b\,\langle \omega \rangle \,,  \label{Apresent}
\end{equation}%
\begin{equation}
h,\,h_{t}\in L_{\infty }(0,T,\,W_{q}^{2}(\Omega ))\hookrightarrow L_{\infty
}(0,T,\,C^{1+\alpha }(\overline{\Omega }))\quad \quad \mbox{
for }\alpha =1-\frac{2}{q}  \label{Bpresent}
\end{equation}%
and $h,\, H$ are the solutions of the systems \eqref{eq4sec2}, respectively.

\underline{2st step.} In this step we prove that the pair $\{\omega ,\mathbf{%
v}\}$ satisfies (\ref{eq7}). Let us define some 
distance functions on $\Gamma $.

\begin{definition}
\label{def2sec copy(1)} Let $d(\mathbf{x},Q):=inf_{\mathbf{y}\in Q}|\mathbf{x%
}-\mathbf{y}|$ be the distance between any given point $\mathbf{x}\in 
\mathbb{R}^{2}$ and any subset $Q\subseteq \mathbb{R}^{2}.$ Let $d(\mathbf{x}%
):=d(\mathbf{x},\mathbb{R}^{2}\backslash \Omega )-d(\mathbf{x},\Omega )$ be
the distance function on $\Gamma $, defined by for any$\ \mathbf{x}\in 
\mathbb{R}^{2}.$

The set of all points of $\overline{\Omega }$, whose distance to $\Gamma $ (
to $\Gamma ^{-}$ and to $\Gamma ^{+}$) is less than $\sigma $, is denoted by 
$U_{\sigma }(\Gamma )$ (respectively, by $U_{\sigma }(\Gamma ^{-})$ and by $%
U_{\sigma }(\Gamma ^{+})$).
\end{definition}

In view of $\Gamma \in C^{2},$ the distance function $d=d(\mathbf{x})$ has
the following properties 
\begin{equation}
d\in C^{2}\quad \quad \mbox{ in }U_{\sigma _{0}}(\Gamma )\quad 
\mbox{ for
some }{\sigma _{0}}>0\quad \quad \mbox{ and }\quad \quad \nabla d=-\mathsf{n}%
\quad \mbox{ on
}\Gamma .  \label{neib}
\end{equation}%
Defining the approximation of the unit function by 
\begin{equation*}
{\mathbf{1}_{\sigma }}(\mathbf{x}):=\left\{ 
\begin{array}{l}
0,\quad \quad \mbox{ if }d(\mathbf{x)\in }[0,\sigma ), \\ 
\frac{d-\sigma }{\sigma },\quad \mbox{ if }d(\mathbf{x)\in }[\sigma ,2\sigma
), \\ 
1,\quad \quad \,\mbox{ if }\mathbf{x}\in \Omega \backslash U_{2\sigma
}(\Gamma ),%
\end{array}%
\right.
\end{equation*}
multiplying the equation of (\ref{eq3sec2}) by $\eta _{\sigma }:={\mathbf{1}%
_{\sigma }}\,\psi $ with $\psi $ being a test function and integrating it
over $\Omega _{T},$ we deduce 
\begin{multline*}
0=\Bigg\{\int_{\Omega _{T}}[b\,\omega _{{\nu }}(\psi _{t}+(\mathbf{v}_{{\nu }%
}\cdot \nabla )\psi )]{\mathbf{1}_{\sigma }-}\left[ \varkappa \omega _{\nu
}+\left( \mathbf{v}_{\nu }\cdot \nabla ^{\bot }\right) \left( \frac{A}{b}%
\right) -\mathrm{rot}\left( \frac{\mathbf{G}}{b}\right) \right] \,\eta
_{\sigma } \\
+\nu \,\omega _{{\nu }}\bigtriangleup \eta _{\sigma }\,d\mathbf{x}%
dt+\int_{\Omega }b\,\omega _{0}(\mathbf{x})\,\eta _{\sigma }(\mathbf{x},0)\,d%
\mathbf{x}\,\Bigg\} \\
+\frac{1}{{\sigma }}\int_{0}^{T}\int_{[\sigma <d<{2\sigma }]}b\,\omega _{{%
\nu }}\,(\mathbf{v}_{{\nu }}\cdot \nabla )d\,\psi \,d\mathbf{x}dt=I^{{\nu }%
,\sigma }+J^{{\nu },\sigma }.
\end{multline*}%
Using (\ref{eq31sec3}), (\ref{V}) and ${\mathbf{1}_{\sigma }}%
\mathop{\longrightarrow}\limits_{{\sigma }\rightarrow 0}1$ in $\Omega _{T}$
and $\Omega $, we have 
\begin{multline*}
\lim_{{\sigma }\rightarrow 0}\,\,\left( \,\lim_{{\nu }\rightarrow 0}I^{{\nu }%
,\sigma }\right) =\int_{\Omega _{T}}\,b\,\omega (\psi _{t}+(\mathbf{v}\cdot
\nabla )\psi ) \\
-\left[ \varkappa \omega +\left( \mathbf{v}\cdot \nabla ^{\bot }\right)
\left( \frac{A}{b}\right) -\mathrm{rot}\left( \frac{\mathbf{G}}{b}\right) %
\right] d\mathbf{x}dt\,+\int_{\Omega }\,b\,\omega _{0}\psi (\mathbf{x},0)\,d%
\mathbf{x}.
\end{multline*}
In view of (\ref{h1}), (\ref{h2}) and a well-known 
\textit{extension} result (see, for instance, p.43, Theorem 3.3, \cite{galdi}%
), there exists an extension $\breve{\omega}_{{\nu }}$ of the boundary
condition $\omega _{\Gamma }(\mathbf{v}_{{\nu }})$ and the initial condition 
$\omega _{0}$ into the domain $\Omega _{T}$, such that 
\begin{equation}
\begin{cases}
\breve{\omega}_{{\nu }}\in W_{q}^{1,1}(\Omega _{T}),\quad \quad & \breve{%
\omega}_{{\nu }}\in L_{\infty }(0,T,\ W_{q}^{2}(\Omega )),\quad \quad \text{
such that } \\ 
||\breve{\omega}_{{\nu }}||_{W_{q}^{1,1}(\Omega _{T})}\leqslant C,\quad \quad
& ||\breve{\omega}_{{\nu }}||_{L_{\infty }(0,T,\ W_{q}^{2}(\Omega
))}\leqslant C\quad \quad \text{ and } \\ 
\breve{\omega}_{{\nu }}\big|_{t=0}=\omega _{0},\quad \quad & \breve{\omega}_{%
{\nu }}\big|_{\Gamma _{T}}=\omega _{\Gamma }(\mathbf{v}_{\nu }),%
\end{cases}
\label{W3}
\end{equation}%
the constants $C=C(\theta ,q)$ do not depend on $\nu $. 
Now, considering  this extension we can  write
\begin{multline*}
J^{\nu ,\sigma }=\frac{1}{\sigma }\int_{0}^{T}\int_{[\sigma <d<2\sigma
]}\,b\,(\mathbf{v}_{{\nu }}\cdot \nabla )d\,(\omega _{\nu }-\breve{\omega}%
_{\nu })\,\psi \,d\mathbf{x}dt \\
+\frac{1}{\sigma }\int_{0}^{T}\int_{[\sigma <d<2\sigma ]}\,b\,(\mathbf{v}_{{%
\nu }}\cdot \nabla )d\,\breve{\omega}_{\nu }\,\psi \,d\mathbf{x}%
dt=J_{1}^{\nu ,\sigma }+J_{2}^{\nu ,\sigma }.
\end{multline*}

  As it was done in 
\cite{C&A} we can show that 
\begin{equation}
\lim_{{\sigma }\rightarrow 0}\,\,\left( \,\overline{\lim_{{\nu }\rightarrow
0}}\,\,\frac{1}{{\sigma }}\int_{0}^{T}\int_{[\sigma <d<{2\sigma }]}b\,\left( 
\mathbf{v}_{\nu }\cdot \nabla \right) d\,|\omega _{{\nu }}-\breve{\omega}_{{%
\nu }}|^{q}\,\,\psi \,d\mathbf{x}dt\,\right) =0,  \label{eq35sec3}
\end{equation}%
for any \underline{positive} test function $\psi $ and any \underline{fixed} 
$q\in (1,\infty ).$ Here $d(\mathbf{x})$ is the distance function,
introduced in  definition \ref{def2sec copy(1)}.

From ({\ref{eq35sec3})} and the boundness of $\omega _{\nu },\,\breve{\omega}%
_{\nu }$ in $L_{\infty }(\Omega _{T})$, independently of $\nu $, we obtain 
\begin{equation*}
\lim_{{\sigma }\rightarrow 0}\,\,\left( \,\overline{\lim_{{\nu }\rightarrow
0}}\,|J_{1}^{{\nu },\sigma }|\right) =0.
\end{equation*}%
By (\ref{h1}), (\ref{h2}) the set of functions $(\mathbf{v}_{\nu }\cdot
\nabla )\,d$ is uniformly continuous on $\overline{\Omega }_{T}$,
independently of ${\nu }$ and the trace of $(\mathbf{v}_{\nu }\cdot \nabla
)\,d$ on $\Gamma _{T}$ is equal to $-a$ for every point $(\mathbf{x},t)\in
\Gamma _{T}.$ By (\ref{W3}) the function $\breve{\omega}_{{\nu }}$ has the
trace $\omega _{{\Gamma }}(\mathbf{v}_{\nu })$\ on the boundary $\Gamma
_{T}^{-}$ \ and, in view of (\ref{w2})-(\ref{h2}) and (\ref{ell-90}),
applied for some $q>3$, we have the convergence 
\begin{equation*}
\omega _{{\Gamma }}(\mathbf{v}_{\nu })\underset{{\nu }\rightarrow 0}{%
\longrightarrow }\omega _{{\Gamma }}(\mathbf{v})\quad \quad \mbox{ in }\quad
C^{\alpha ,\alpha }(\Gamma _{T})\quad \mbox{ with  }\quad \alpha <1-\frac{3}{%
q},
\end{equation*}%
that gives 
\begin{equation*}
\lim_{{\sigma }\rightarrow 0}\,\,(\lim_{{\nu }\rightarrow 0}\,\,J_{2}^{{\
\nu },\sigma })=-\int_{0}^{T}\int_{\Gamma _{T}^{-}}\omega _{{\Gamma }}(%
\mathbf{v})\,a\,\psi \,d\mathbf{x}dt.
\end{equation*}%
Therefore the pair $\{\omega ,\,\mathbf{v}\}$ satisfies the equation (\ref%
{eq7}). \bigskip $\hfill \blacksquare $

\vspace{1pt}
Let us introduce the  sets
 $\Omega ^{\theta }:=\left\{ \mathbf{x}\in \mathbb{R}^{2}:\ d(\mathbf{x}%
,\Omega \mathbf{)}<\theta \right\}  \,$ and $\,\Omega
_{T}^{\theta }:=\Omega ^{\theta }\times \left[ -\theta ,T+\theta \right] .$
In view of the \textit{extension} result (p.43, Theorem 3.3, \cite{galdi})
and (\ref{ap2}), since $\mathbf{v}\in W_{q}^{1,1}(\Omega _{T}),$ there
exists an {\ extension} $\breve{\omega}$ of the boundary conditions $\omega
_{\Gamma }(\mathbf{v}):=\gamma \,\mathbf{v}\cdot {\mathsf{s}}+g$ and the
initial conditions $\omega _{0}$ into the domain $\Omega _{T}^{\theta }$,
such that for any $q\in (1,\infty )$ 
\begin{equation*}
\left\{ 
\begin{array}{ll}
\breve{\omega}\in W_{q}^{1,1}(\Omega _{T}^{\theta }),\quad \quad 
\mbox{ with
} & ||\breve{\omega}||_{W_{q}^{1,1}(\Omega _{T})}\leqslant C\quad \quad 
\mbox{
and } \\ 
&  \\ 
\breve{\omega}\big|_{t=0}=\omega _{0}, & \breve{\omega}\big|_{\Gamma
_{T}}=\omega _{\Gamma }(\mathbf{v}),%
\end{array}%
\right.
\end{equation*}%
where the  constant $C$ depends on $\theta ,q$ and $||\mathbf{v||}%
_{W_{p}^{1,1}(\Omega _{T})}.$ By the same approach as in \cite{C&A} we can
deduce

\begin{align}
\lim_{{\sigma }\rightarrow 0}\,\,\left( \,\,\frac{1}{{\sigma }}%
\int_{0}^{T}\int_{[\sigma <d<{2\sigma }]}b\,(\mathbf{v}\cdot \nabla
)\,d\,|\omega -\breve{\omega}|^{q}\,\psi \,d\mathbf{x}dt\,\right) & =0,
\label{bo} \\
\lim_{{\sigma }\rightarrow 0}\,\,\left( \frac{1}{{\sigma }}\int_{0}^{{\sigma 
}}\int_{\Omega }b\,|\omega -\breve{\omega}|^{q}\,\,\psi \,d\mathbf{x}%
dt\right) & =0  \label{to}
\end{align}%
for any \underline{positive} test function $\psi $ and any $q\in (1,\infty
). $
 With the help of \eqref{bo} \eqref{to} and  applying the methods in \cite{C&A}, we obtain the
following Gronwall's type inequality$:$

\begin{teo}
\label{lem6sec4234} Let $q\in (1,\infty )$ be given. For all $t\in \lbrack
0,T]$, we have 
\begin{align}
\int_{\Omega }b\,|\omega (\mathbf{x},t)|^{q}\,d\mathbf{x}& -\int_{\Omega
}b\,|\omega _{0}|^{q}\,d\mathbf{x}\leqslant q\int_{0}^{t}\int_{\Omega
}(|A|+|\varkappa |)\,|\omega |^{q}\,d\mathbf{x}\,d\tau  \notag \\
& +q\int_{0}^{t}\left\vert \mathrm{rot}\left( \frac{\mathbf{G}}{b}\right)
\right\vert |\omega |^{q-1}\ d\mathbf{x}d\tau +\int_{0}^{t}\left\vert (%
\mathbf{v}\cdot \nabla ^{\perp })\left( \frac{A}{{b}}\right) \right\vert
\,|\omega |^{q-1}\,d\mathbf{x}\,d\tau  \notag \\
& +\int_{0}^{t}\int_{\Gamma ^{-}}b\,a\,|\omega _{\Gamma }(\mathbf{v})|^{q}\,d%
\mathbf{x}d\tau .  \label{A}
\end{align}
\end{teo}

\bigskip

\section{Limit transition on $\protect\theta $}

\setcounter{equation}{0}

\label{secbbb}

In the previous sections we have constructed the solution $\{\omega ,h,%
\mathbf{v}\}$ for the system (\ref{eq7}), (\ref{eq4sec2}), depending on ${%
\theta }$. In the sequel instead of $\omega ,\,h,\,\mathbf{v}$ and the
approximated data $a,\,\gamma ,\,g,\,\omega _{0},\,\varkappa ,\,A,\,\mathbf{G%
}$ we shall write $\omega _{\theta },$ $h_{\theta },\mathbf{v_{\theta }}\ $
and $a^{\theta },\,\gamma ^{\theta },\,g^{\theta },\,\omega _{0}^{\theta
},\,\varkappa ^{\theta },\,\,A^{\theta },\,\mathbf{G}^{\theta }$,
respectively.

\textit{In this subsection all constants } $C$ \textit{\ }\underline{\textit{%
do not depend on } $\theta $}.

Let us formulate a  Gronwall's type lemma. This lemma can be proved by the
standard method.

\begin{lemma}
\label{lem4.0} Let $D(t),B(t)\in L_{1}(0,T)$ be given {\ non-negative }
functions. Let $y(t)$ be a non-negative function for $t\in \lbrack 0,T]$ and 
$y(t)=0$ for $\forall t<0$, satisfying 
\begin{equation}
y(t)\leqslant y_{0}+\int_{0}^{t}\left[ D(\tau )\cdot (u(\tau )+y(\tau
))+B(\tau )\right] d\tau  \label{gronwall}
\end{equation}%
%
%
%
%
%
%
%
%
%
%
%
%
%
%
%
%
%
%
%
%
%
%
%
%
%
%
%
%
%
%
%
%
%
%
%
%
%
%
%
%
%
%
%
%
%
%
%
%
%
%
%
%
%
%
%
%
%
%
with $u(t):=\frac{1}{\theta }\int_{t-\theta }^{t}y(\tau )\,d\tau $. Then
there exists $\theta _{0}>0$, such that for any fixed $\theta \in (0,$ $%
\theta _{0})$ we have 
\begin{equation*}
y(t)\leqslant 2\exp \left( \int_{0}^{t}D(\tau )d\tau \right) \left[
y_{0}+\int_{0}^{t}B(r)\cdot \exp \left( -\int_{0}^{r}D(\tau )d\tau \right) dr%
\right] ,\quad \quad \forall t\in \lbrack 0,T].
\end{equation*}
\end{lemma}

\bigskip

Below  the limit transitions on the regularization parameter $\theta 
$ will be done first for  $p\in (1,\infty)$. The case $p=\infty$ will
be considered just at the end of the article. Combining Theorem \ref%
{lem6sec4234} and\ Lemma \ref{lem4.0}, we show the following result:

\begin{lemma}
\label{teo4sec2 copy(1)} The pair $\{\omega _{\theta },h_{\theta }\}$
satisfies the estimates 
\begin{align}
||\omega _{\theta }||_{L_{\infty }(0,T,\,L_{p}(\Omega ))}& \leqslant C,
\label{t1} \\
||h_{\theta }||_{L_{\infty }(0,T,\,W_{p}^{2}(\Omega ))}& \leqslant C,
\label{t2}
\end{align}%
for a.e. $t\in (0,T)$ and any $|\Delta |<\min \{t,T-t\}$.
\end{lemma}

\textbf{Proof.} Let us denote by $\,y(t):=||\omega_\theta(\cdot,t)||^p_{L_p(%
\Omega)} \,$ and $\, u(t):= \frac{1}{\theta }\int_{t-\theta}^{t}y(\tau)\,d%
\tau $. Now, we estimate different terms of the right member in \eqref{A}.

$1^{st}$ term: By \eqref{eq00sec1}-\eqref{eqlinha} and \eqref{ell-90}, we
have $A\in L_{2}(0,T,W_{\widetilde{p}}^{1}(\Omega ))\hookrightarrow
L_{2}(0,T,C(\bar{\Omega})),$ since $\widetilde{p}>2\,.$ Taking into account %
\eqref{ap2}, this implies 
\begin{equation}\label{1term}
\int_{0}^{t}\int_{\Omega }(|A^{\theta }|+|\varkappa ^{\theta }|)\,|\omega
_{\theta }|^{p}dx\,d\tau \leqslant \int_{0}^{t}f_{1}(\tau )\,\,y(\tau
)\,d\tau ,
\end{equation}%
where $f_{1}:=||A||_{C(\bar{\Omega})}+||\varkappa ||_{L_{\infty }(\Omega
)}\in L_{1 }(0,T)$.

$2^{nd}$ term: By \eqref{eq00sec1} and \eqref{ap2}, we see that 
\begin{align}
\int_0^t \int_\Omega |\mathrm{rot}\left(\frac{\mathbf{G^\theta }}{b}%
\right)|\,|\omega_\theta |^{p-1}\, dx\,d\tau\leqslant\int_0^t f_2 (\tau)
\,\, y(\tau)^{1-1/p} \, d\tau  \notag \\
\leqslant C + \int_0^t f_2 (\tau) \,\, y(\tau)\, d\tau  \label{ast}
\end{align}
where $f_2 :=||\mathrm{rot}\left(\frac{\mathbf{G}}{b}\right)\,||_{L_p(%
\Omega)} \in L_1 (0,T)$. Here we used 
$y^{1-1/p} \leqslant y +1$.

$3^{d}$ term: We have that 
\begin{equation}  \label{I}
\begin{split}
I&:=\int_0^t\int_\Omega \left|(\mathbf{v_\theta }\cdot\nabla)\left(\frac{%
A^\theta }{b}\right)\right|\, | \omega_\theta |^{p-1}dx\,d\,\tau \\
&\leqslant \int_0^t f_3 (\tau) \, \|\mathbf{v}(\cdot,\tau)\|_{L_{p_2}(%
\Omega)} \,\, y(\tau)^{1-\frac{1}{p}}\,d\,\tau ,
\end{split}%
\end{equation}
where $f_3 :=\|\nabla\left(\frac{A}{b}\right)\|_{L_{p_1}(\Omega)}\in
L_2 (0,T)$. Here we used the H\" older inequality for $p_1:=\widetilde{p%
}$ (defined in \eqref{eqlinha}), 
\begin{equation}  \label{p2}
p_2:= 
\begin{cases}
\infty, \quad \quad \, \, \mbox{ if }\quad p\in(2, \infty), \\ 
2+\frac{4}{\epsilon}, \quad \, \mbox{ if }\quad p=2, \\ 
\frac{p}{2-p}, \quad \quad \mbox{ if }\quad p\in ( 1, 2)%
\end{cases}%
\end{equation}
and $p_3:=\frac{p}{p-1}$, which satisfy the identity $\frac{1}{p_1}+ \frac{1%
}{p_2}+\frac{1}{p_3}=1$.

In view of \eqref{eq00sec1}, \eqref{ap2}, \eqref{eq4sec2}, \eqref{VELOC} and %
\eqref{V}, \eqref{Apresent}, we have for a.e. $\tau \in (0,T)$ 
\begin{equation}  \label{for v}
\begin{split}
\| \mathbf{v}_\theta (\cdot,\tau) \|_{W_p^1 (\Omega)} \leqslant C ( \|b \,
\langle \omega_\theta \rangle (\cdot , \tau ) \|_{L_p (\Omega)} +f_4 (\tau)
) &\leqslant C\,( u^{\frac{1}{p}}(\tau ) +f_4 (\tau ) ) \,,
\end{split}%
\end{equation}
where $f_4 :=\|A \|_{L_{p}(\Omega)} +\|a \|_{W_p^{1-\frac{1}{p}}(\Gamma)}\in
L_2 (0,T) $.

By the last inequality and the embedding $W^1_{p}(\Omega) \hookrightarrow
L_{p_2}(\Omega) $, following from \eqref{ell-90}, we derive 
\begin{equation*}  \label{v11}
\begin{split}
\|\mathbf{v}_\theta (\cdot,\tau)\|_{L_{p_2}(\Omega)}&\leqslant C\biggl( u^{%
\frac{1}{p}}(\tau ) +f_4 (\tau ) \biggr) \quad \quad \text{ for a.e. } \tau
\in (0,T).
\end{split}%
\end{equation*}
Hence we conclude 
\begin{equation}  \label{III}
\begin{split}
I \leqslant C\int_0^t f_5 (\tau ) ( u(\tau ) +y(\tau ) )\, d\tau +C
\end{split}%
\end{equation}
for  $f_5 \in L_1 (0,T)$. Here we applied the inequality $y^{1-1/p} \leqslant y +1$
 and 
$u^{1/p} \, y^{1-1/p} \leqslant C_p (u +y) $ following from \eqref{eqest}.

$4^{th}$ term: By \eqref{ell-9}, we have $W^1_{p}(\Gamma) \hookrightarrow
C(\Gamma) $, therefore, using \eqref{eq00sec1}-\eqref{eqlinha} and %
\eqref{ap2} we can estimate 
\begin{align*}
J:&=\int_0^t\int_{\Gamma^-} a^\theta\,\left| \omega_\Gamma(\mathbf{v}%
_\theta) \right|^p dx\,d\,\tau  \notag \\
&\leqslant C  \int_0^t\biggl(%
\int_{\Gamma^-}  \|a\|_{C(\Gamma ) }\, \left|\gamma^\theta \right|^p \left|\mathbf{v}_\theta
\right|^p \, dx+\|g^\theta\|_{L_p(\Gamma^-)}^p \biggr)d\,\tau  \notag \\
&\leqslant C\int_0^t \|a\|_{C(\Gamma ) }\, \|\gamma \|^p_{L_{\widetilde p}(\Gamma)} \,\cdot \| \mathbf{v%
}_\theta \|^p_{L_{p_2}(\Gamma)}\,d\tau+C,  
\end{align*}
where $\widetilde{p}$ and ${p_2}$ are defined by \eqref{eqlinha} and \eqref{p2},
respectively. By \eqref{ell-9} we have the embedding $W^1_{p}(\Omega) \hookrightarrow
L_{p_2}(\Gamma) $, that implies from \eqref{for v} and %
\eqref{eq00sec1}-\eqref{eqlinha} 
\begin{equation}\label{JJJ}
J\leqslant C\int_0^t f_6(\tau)\,u(\tau)\,d\,\tau+C
\end{equation}
with $f_6:=\|a\|_{C(\Gamma ) }\, \|\gamma \|_{L_{\widetilde p}(\Gamma)}^p\in L_1 (0,T) .$

Combining all estimates derived for different terms \eqref{1term},\eqref{ast}, \eqref{III}, \eqref{JJJ}  and applying Lemma \ref{lem4.0}, we derive the 
desired estimate \eqref{t1}. Estimate \eqref{t2} follows from \eqref{ell-1}.
 $\hfill \;\blacksquare $

\begin{lemma}
\label{7}There exists a constant $C>0$, independent of $\theta $, such that: 
\begin{equation}
||\partial _{t}\omega _{\theta }||_{L_{\infty }(0,T;\,\,W_{p}^{-2}(\Omega
))}\leq C.  \label{22e}
\end{equation}
\end{lemma}

\vspace{1pt}\textbf{Proof.} Let us choose in (\ref{eq7}) the test function $%
\psi (\mathbf{x},t):=\phi (\mathbf{x})\,\varphi (t)$, such that $\phi (%
\mathbf{x})\in C^{2}(\Omega )$ and $\varphi (t)\in C^1 (0,T)$\ \ and $%
\mbox{
supp}\,(\psi )\subset \Omega _{T}.$ Then from \eqref{eq7}, using (\ref{V}), (%
\ref{eq66sec2}), the definition (\ref{WW}) and the symmetry of the kernel $%
K_{1}$ on the variables $\mathbf{x},\,\mathbf{y}$, we obtain 
\begin{eqnarray}
\int_{0}^{T}\biggl(\int_{\Omega }\omega _{\theta }\,\phi \,d\mathbf{x}\biggl)%
\varphi _{t}\,dt &=&-\int_{0}^{T}\biggl(\int_{\Omega }\,\omega _{\theta }%
\mathbf{v}_{\theta }\cdot \nabla \phi \,d\mathbf{x}\biggr)\varphi
\,dt=-F[\omega _{\theta },K_{\psi }]  \notag \\
&-&\int_{0}^{T}\varphi \left( \int_{\Omega }b\omega _{\theta }\nabla
H_{\theta }\cdot \nabla \phi \,d\mathbf{x\varphi }\right) dt\mathbf{:=}F+I.
\label{22ee}
\end{eqnarray}%
Since $\Gamma $ is $C^{2}$-smooth, the kernel $K_{1}$ satisfies the
inequality $\left\vert \nabla _{\mathbf{x}}K_{1}(\mathbf{x},\mathbf{y}%
)\right\vert \leqslant C|\mathbf{x}-\mathbf{y}|^{-1}$ for any $\mathbf{x},%
\mathbf{y}\in \Omega .$ Therefore from (\ref{t1}), $b\in C^{1}(\overline{%
\Omega })$ and Theorem 1 of \cite{hajlasz}, we have 
\begin{eqnarray} \label{II2}
|F| &\leqslant &C\Vert \omega _{\theta }\Vert _{L_{p}(\Omega )}^{2}\left(
\int_{\Omega }\int_{\Omega }\left( \frac{|b(\mathbf{y})\,\nabla \phi (%
\mathbf{x})-b(\mathbf{x})\,\nabla \phi (\mathbf{y})|}{|x-y|}\right)
^{p^{\ast }}d\mathbf{x}\,d\mathbf{y}\right) ^{1/p^{\ast }}||\varphi
||_{L_{1}(0,T)}  \notag   \\
&\leqslant &C\Vert \phi \Vert _{W_{p^{\ast }}^{2}(\Omega )}||\varphi
||_{L_{1}(0,T)}\quad \;\text{with}\;\quad \frac{1}{p^{\ast }}+\frac{1}{p}=1.
\end{eqnarray}%
Considering three diferent cases $p>2,$ $p=2$ and $1<p<2,$ the term $I$ is
estimated as 
\begin{equation}
|I|\leqslant \Vert \omega _{\theta }\Vert _{L_{p}(\Omega )}\Vert \nabla
H_{\theta }\cdot \nabla \phi \Vert _{L_{p^{\ast }}(\Omega )}^{2}||\varphi
||_{L_{1}(0,T)}\leqslant C\Vert \phi \Vert _{W_{p^{\ast }}^{2}(\Omega
)}||\varphi ||_{L_{1}(0,T)}.  \label{II3}
\end{equation}%
Let us note that for $1<p<2,$ we have $A(\cdot ,t)\in W_{\widetilde{p}%
}^{1}(\Omega )\hookrightarrow C^{\delta }(\overline{\Omega })$ and $a(\cdot
,t)\in W_{p}^{1}(\Gamma )\hookrightarrow C^{\delta }(\Gamma )$ for some $%
\delta >0,$ then from the theory of potentials (see \cite{vlad}) we derive
that $\Vert \nabla H_{\theta }(\cdot ,t)\Vert _{C(\overline{\Omega }%
)}\leqslant C,$ independently of $\theta $ and $t.$

The  relations \eqref{22ee}-\eqref{II3} imply 
\begin{equation*}
\Bigl|\int_{0}^{T}\Bigl(\int\limits_{\Omega }\omega _{\theta }\,\phi \;d%
\mathbf{x}\Bigl)\varphi _{t}dt\Bigl|\leqslant C\,\Vert \phi \Vert
_{W_{p^{\ast }}^{2}(\Omega )}||\varphi ||_{L_{1}(0,T)},
\end{equation*}%
which is equivalent to (\ref{22e}). \bigskip $\hfill \;\blacksquare $

From (\ref{t1})-(\ref{t2}), the representation (\ref{Apresent}) and the
approximation convergence (\ref{ap3}), we conclude that there exists a
subsequence of $\{\omega _{\theta },h_{\theta },\mathbf{v}_{\theta }\}$,
such that 
\begin{align}
h_{\theta }\rightharpoonup h\quad \quad & \mbox{ weakly}-\ast \mbox{ in }%
L_{\infty }(0,T,\,W_{p}^{2}(\Omega )),  \notag \\
\omega _{\theta }\rightharpoonup \omega \quad \quad & \mbox{ weakly}-\ast 
\mbox{ in
}L_{\infty }(0,T,L_{p}(\Omega )),  \label{transition} \\
\mathbf{v}_{\theta }\rightharpoonup \mathbf{v}\quad \quad & \mbox{ weakly}%
-\ast \mbox{ in
}L_{\infty }(0,T,W_{p}^{1}(\Omega ))  \notag
\end{align}%
and $\,\mathbf{v}\,$ satisfies the relation \eqref{V}.

The following lemma \ plays the main role in the  argument for 
the limit transition.

\begin{lemma}
\label{lem33}For arbitrary test function $\psi ,$ we have 
\begin{equation}
F[\omega _{\theta },K_{\psi }]\rightarrow F[\omega ,K_{\psi }]\;\;\;\,\;\;%
\mbox{ as  }\;\;\;\,\;\;\theta \rightarrow 0.  \label{WWZ}
\end{equation}
\end{lemma}

\noindent

\textbf{Proof.} Let ${\mathbf{\rho }}(s)\in C_{0}^{\infty }(\mathbb{R})$
with $\mathbf{\rho } (s):= 1,$ if $|s|\leqslant 1$ and $0$, if $|s|>2$. We
introduce the functions ${\mathbf{\rho }_{\sigma }}(\mathbf{x}):={\mathbf{%
\rho }}(|\mathbf{x}|/\sigma )$ and ${\mathbf{\rho }}_{\sigma }^{\Gamma }(%
\mathbf{x}):={\mathbf{\rho }}(d(\mathbf{x})/\sigma )$ for $\sigma >0$.

We can show that 
\begin{equation}
\sup_{\theta }|F[\omega _{\theta },z_{\sigma }]| \underset{{\sigma
\rightarrow 0}}{\longrightarrow } 0  \label{Wlim}
\end{equation}
for each $z_{\sigma }(\mathbf{x},\mathbf{y}):=\mathbf{\rho }_{\sigma
}^{\Gamma }(\mathbf{x})\ K_{\psi }(\mathbf{x},\mathbf{y})$, $z_{\sigma }(%
\mathbf{x},\mathbf{y}):=(1-\mathbf{\rho }_{\sigma }^{\Gamma }(\mathbf{x}))%
\mathbf{\rho }_{\sigma }^{\Gamma }(\mathbf{y})\ K_{\psi }(\mathbf{x},\mathbf{%
y})$ and $z_{\sigma }(\mathbf{x},\mathbf{y}):=(1-\mathbf{\rho }_{\sigma
}^{\Gamma }(\mathbf{x}))\mathbf{\rho }_{\sigma }(\mathbf{x-y})\ K_{\psi }(%
\mathbf{x},\mathbf{y})$. For instance if $z_{\sigma }(\mathbf{x},\mathbf{y})=\mathbf{\rho }_{\sigma
}^{\Gamma }(\mathbf{x})\ K_{\psi }(\mathbf{x},\mathbf{y})$ as in the estimate \eqref{II2} we have
\begin{eqnarray}
|F[\omega _{\theta },z_{\sigma }] | \leqslant C \left(
\int_{\Omega } \int_{\Omega }\left( \mathbf{\rho }_{\sigma
}^{\Gamma }(\mathbf{x})  \frac{|b(\mathbf{y})\,\nabla \phi (%
\mathbf{x})-b(\mathbf{x})\,\nabla \phi (\mathbf{y})|}{|x-y|}\right)
^{p^{\ast }}d\mathbf{x}\,d\mathbf{y}\right) ^{1/p^{\ast}}||\varphi
||_{L_{1}(0,T)},  \notag   
\end{eqnarray}%
therefore, thank to the continuity property of integral, we obtain \eqref{Wlim}.

The compact embedding of $L_{p}(\Omega )$ in $%
W_{p}^{-2}(\Omega )$ and (\ref{22e}) imply 
\begin{equation}
\omega _{\theta }\underset{{\theta 
\rightarrow 0}}{\longrightarrow} \omega \quad \mbox{strongly in}\quad
L_{2}(0,T;\;W_{p}^{-2}(\Omega ))
\end{equation}%
by results obtained in \cite{aubin}, \cite{simon}. Since $z_{\sigma }:=(1-%
\mathbf{\rho }_{\sigma }^{\Gamma }(\mathbf{x}))(1-\mathbf{\rho }_{\sigma
}^{\Gamma }(\mathbf{y})-\mathbf{\rho }_{\sigma }(\mathbf{x}-\mathbf{y}%
))K_{\psi }(\mathbf{x},\mathbf{y})\in C^{2}(\overline{\Omega }\times 
\overline{\Omega })$ is a smooth function with a compact support in $\Omega
\times \Omega $, we have 
\begin{equation*}
F[\omega _{\theta },z_{\sigma }]\underset{{\theta 
\rightarrow 0}}{\longrightarrow}  F[\omega ,z_{\sigma }],
\end{equation*}%
that, jointly with \eqref{Wlim}, implies \eqref{WWZ}$.$ \bigskip $\hfill
\;\blacksquare $

Using (\ref{transition})-(\ref{WWZ}), the limit transition on $\theta
\rightarrow 0$ \ in the equations (\ref{eq7}) and (\ref{eq4sec2}), written
for $\{\omega _{\theta },h_{\theta },\mathbf{v}_{\theta }\}$\vspace{1pt},
implies that the limit triple $\{\omega ,\,h,\,\mathbf{v}\} $\ satisfies (%
\ref{eq8})-(\ref{eq77}).

\bigskip

We investigated the case $p\in (1,\infty )$. Let us study the case $%
p=\infty .$  We consider approximated  data $\gamma ^{\theta },$ $g^{\theta },$
$A^{\theta }$, $\omega _{0}^{\theta }$, $a^{\theta }$\ satisfying the
conditions (\ref{ap2}) for any $q\in (2,\infty ]$ 
and (\ref{ap3})  for any $q\in (2,\infty )$ (instead of $p$). By the above
construction, there exists an approximate solution $\{\omega _{\theta
},h_{\theta },\mathbf{v}_{\theta }\},$ satisfying\ the estimates of Lemma %
\ref{teo4sec2 copy(1)} and Lemma \ref{7}. These estimates hold for any $%
q\in (2,\infty )$ (instead of  $p$), but for constants $C$, depending
on $q$.

Let us return to deducing the  estimate (\ref{t1}). From the estimate (%
\ref{VELOC}), written for some fixed $p=\overline{q}\in (2,\infty )$ and the
embedding $W_{\overline{q}}^{1}(\Omega )\hookrightarrow C(\overline{\Omega }%
),$ we have 
\begin{equation}
\Vert \mathbf{v}_{\theta }\Vert _{L_{\infty }(0,T,C(\overline{\Omega }%
))}\leqslant C\Vert \mathbf{v}_{\theta }\Vert _{L_{\infty }(0,T,W_{\overline{%
q}}^{1}(\Omega ))}\leqslant C(\overline{q}),  \label{inf1}
\end{equation}%
where $C(\overline{q})$ depends only on $\overline{q}.$ Therefore estimating
the $1^{st}-3^{d}$ terms\ of (\ref{A}) by  the same way as in Lemma \ref%
{teo4sec2 copy(1)} and  the $4^{th}$ term using (\ref{t1}), we deduce that 
$\Vert \omega _{\theta }(\cdot, t )\Vert _{L_{q}(\Omega )} \leqslant \Vert
\omega _{0}\Vert _{L_{q}(\Omega )}+C(\overline{q})$ for a.e. $t\in (0,T)$.
Hence the limit $q\rightarrow \infty $ gives 
\begin{equation*}
\Vert \omega _{\theta }\Vert _{L_{\infty }(\Omega_T)}\leqslant \Vert \omega
_{0}\Vert _{L_{\infty }(\Omega )}+C(\overline{q}).
\end{equation*}%
This  estimate together with (\ref{inf1}) allows to pass to the limit in $%
\theta \rightarrow 0$ in (\ref{eq7}) and (\ref{eq4sec2}), written for $%
\{\omega _{\theta },h_{\theta },\mathbf{v}_{\theta }\}$, which proves Theorem %
\ref{teo2sec1} for the case $p=\infty $. \bigskip $\hfill \;\blacksquare $

\section*{Acknowledgement}

N.V. Chemetov thanks the support from FCT and FEDER through the Project
\linebreak POCTI/209/2003 of Centro de Matem\'{a}tica e Aplica\c{c}\~{o}es
Fundamentais da Universidade de Lisboa (CMAF/UL). The work of F. Cipriano is
partially supported by the Portuguese research project POCTI/MAT/55977/2004.

\end{document}